\documentclass[12pt]{article}
\usepackage{amsmath,amsfonts,amsthm,amssymb,eucal,amsxtra}
\newcommand{\B}{\mathcal B}
\newcommand{\K}{\widetilde{K}}
\newcommand{\M}{\mathbf M}
\newcommand{\AAA}{\mathfrak A}
\newcommand{\N}{\mathbf N}
\newcommand{\A}{\widetilde{\mathfrak A}}
\newcommand{\LL}{\mathbf L}
\newcommand{\RR}{\mathbf R}
\newcommand{\F}{\mathcal F_p}
\begin{document}
\newtheorem*{lem}{Lemma}
\newtheorem*{teo}{Theorem}
\newtheorem*{prop}{Proposition}
\pagestyle{plain}
\title{Non-Archimedean Group Algebras with Baer Reductions}
\author{Anatoly N. Kochubei\\
\footnotesize Institute of Mathematics,\\
\footnotesize National Academy of Sciences of Ukraine,\\
\footnotesize Tereshchenkivska 3, Kiev, 01601 Ukraine
\\ \footnotesize E-mail: \ kochubei@i.com.ua}
\date{}
\maketitle

\bigskip
\begin{abstract}
Within the concept of a non-Archimedean operator algebra with the Baer reduction (A. N. Kochubei, On some classes of non-Archimedean operator algebras, Contemporary Math. {\bf 596} (2013), 133--148), we consider algebras of operators on Banach spaces over non-Archimedean fields generated by regular representations of discrete groups.
\end{abstract}

\medskip
{\bf MSC 2010}. Primary: 47S10. Secondary: 47L10; 43A99; 16S34.

\bigskip
\section{INTRODUCTION}

While most of traditional branches of analysis have well-developed non-Archimedean analogs, very little is known about non-Archimedean operator algebras. The reason is the absence of natural involutions on algebraically closed non-Archimedean fields and related structures. An involution is vital for the theories of $C^*$- and $W^*$ (von Neumann)-algebras, which is stressed even in the names of these objects. Therefore the methods developed in the theory of operator algebras do not carry over to the non-Archimedean situations.

An alternative approach to the problem was suggested in \cite{K}. Let $K$ be a complete non-Archimedean valued field with a nontrivial valuation, $\K$ be the corresponding residue field (see \cite{BGR,PS,R} for the main notions and results of non-Archimedean analysis). Let us consider a unital algebra $\AAA$ of bounded linear operators on a Banach space $\B$ over $K$ possessing an orthonormal basis (in the non-Archimedean sense). In this case $\B$ is isomorphic to the space $c_0(J,K)$ of sequences $x=[x_i;i\in J]$, $x_i\in K$, which, if the set $J$ is infinite, tend to 0 by the filter of complements to finite subsets of $J$. We may assume that $\B =c_0(J,K)$. The algebra $\AAA$ is endowed with the usual norm of operators. Denote by $\AAA_1$ the closed unit ball in $\AAA$ -- the set of operators with the norm $\le 1$.

$\AAA_1$ is an algebra over the ring $O$ of integers of the field $K$, as well as its ideal $\AAA_0$ consisting of operators of norm $<1$. {\it The reduction} $\A =\AAA_1/\AAA_0$ can be considered as a $\K$-algebra with unit. Note that elements of $\AAA_1$ can be interpreted as finite or infinite matrices with entries from $O$; see \cite{Se}. Then $\A$ consists of matrices with entries from $\K$.

In order to find a non-Archimedean counterpart of the class of von Neumann algebras, we were looking in \cite{K} for a class of $\K$-algebras, for which there is a (purely algebraic) theory parallel to the theory of von Neumann algebras. Such is the class of Baer rings and algebras introduced by Kaplansky \cite{Kap}. A unital ring $R$ is called a Baer ring, if each left (or, equivalently, each right) annihilator in $R$ is generated by an idempotent element. This property, proved by Baer \cite{B} for the ring of all endomorphisms of a vector space of an arbitrary dimension, is valid for some classes of rings including von Neumann algebras.

A Baer ring $R$ is called Abelian, if all its idempotents are central, and Dedekind finite, if $xy=1$ implies $yx=1$. An idempotent $e\in R$ is called Abelian (finite), if the Baer ring $eRe$ is Abelian (respectively, Dedekind finite). If $u$ and $v$ are central idempotents, we write $u\le v$, if $vu=u$. An idempotent $e$ is called faithful, if the smallest of the central idempotents $v$, such that $ve=e$, is equal to 1.

Kaplansky \cite{Kap} introduced the following types of Baer rings. A Baer ring $R$ is of type I, if it has a faithful Abelian idempotent. It is of type II, if it has a faithful finite idempotent, but no nonzero Abelian idempotents, and of type III, if it has no nonzero finite idempotents. These classes are subdivided further into finite and infinite ones. The main result \cite{Kap} of the theory of Baer rings is the unique decomposition of every Baer ring into a finite direct sum of rings of the above types. By the construction, for a Baer algebra, the summands are algebras too.

We call an operator algebra $\AAA$ {\it an algebra with the Baer reduction}, if its reduction $\A$ is a Baer ring. It is known (see \cite{BGR}, Lemma 2.5.1/3) that a finite system of elements of norm 1 in a non-Archimedean normed space is orthonormal if and only if their reductions are linearly independent. Therefore the operator ring $\AAA_1$ with the Baer reduction is an orthogonal sum of rings with reductions of types I, II, and III.

More generally, we can consider the situation where $\A$ is not necessarily a Baer ring, but its maximal right ring of quotients $Q(\A )$ is (see e.g. \cite{F,L} for the definition and properties of this ring). Note that the ring $Q(\A )$ constructed from the algebra $\A$ is actually a $\K$-algebra too (see Lemma 14.15 in \cite{L}), and its decomposition into a direct sum of subalgebras also generates via the monomorphism $\A \to Q(\A )$ a decomposition of $\A$, thus the orthogonal decomposition of $\AAA_1$. In this situation we say that $\AAA$ is {\it an algebra with weak Baer reduction}.

Thus, in both the above situations we have a kind of classification of operator algebras. In order to justify the above definitions, we need meaningful examples with (possibly weak) Baer reductions. In \cite{K}, we developed a version of the crossed product construction for actions of a class of compact groups with good $p$-adic harmonic analysis.

In this paper, we follow the idea (and some techniques) of \cite{MN} and consider analogs of the von Neumann group algebras, operator algebras generated by the right and left regular representations of a discrete group $G$. For this case, if, for example, $K=\mathbb C_p$ (the completion of an algebraic closure of the field $\mathbb Q_p$ of $p$-adic numbers), then $\K =\mathcal F_p$ (an algebraic closure of the finite field $\mathbb F_p$), and $\A =\mathcal F_p[G]$ is the group algebra over $\mathcal F_p$ in the usual sense \cite{P}, and the existing results about the latter produce the desired examples.

The author is grateful to the anonymous referee for helpful comments.

\section{Operator Topologies}

Let $A$ be a bounded linear operator on $\B =c_0(J,K)$. Any element $x\in \B$ can be written as a convergent series
$$
x=\sum\limits x_j\delta_j
$$
where $\delta_j=\left( \delta_{j,i}\right)_{i\in J}$, and $\delta_{j,i}$ is the Kronecker symbol. Then $A$ corresponds to a matrix (in general, infinite) $\left( a_{i,j}\right)$, $a_{i,j}\in K$, such that the double sequence $|a_{i,j}|$ is bounded ($|\cdot |$ is the absolute value in $K$) and, if $J$ is infinite, $a_{i,j}\to 0$ for any fixed $j$ and $i\to \infty$ (by the filter of complements to finite sets); moreover,
$$
\| A\| =\sup\limits_{(i,j)}\left| a_{i,j}\right| ;
$$
see \cite{Se}. As usual \cite{BR}, the norm topology on the set $L(\B )$ of all bounded linear operators on $\B$ is called {\it the uniform operator topology}. The topology defined by the set of all the seminorms $\| A\|_x=\|Ax\|$, $x\in \B$, is called {\it the strong operator topology}.

In the non-Archimedean situation, the strong operator topology can be described in terms of matrix elements of operators.

\medskip
\begin{prop}
Let $(\mathcal L,\succ )$ be a directed set, $A^{(l)}= \left( a_{i,j}^{(l)}\right)$, $l\in \mathcal L$, be a net of operators from $L(\B)$. Then the convergence $A^{(l)}\to 0$ in the strong operator topology is equivalent to the following two properties:

\begin{equation}
\| A^{(l)}\| \le C,\ C>0,\ \text{for all $l\in \mathcal L$};
\end{equation}

\begin{equation}
\sup\limits_i\left| a_{i,j}^{(l)}\right| \to 0\quad \text{for each fixed $j$}.
\end{equation}
\end{prop}

\medskip
{\it Proof}. If $A^{(l)}\to 0$ in the strong operator topology, then $A^{(l)}\delta_j\to 0$ strongly in $c_0(J,K)$ for each $j$. Since $A^{(l)}\delta_j=\left[ a_{i,j}^{(l)},i\in J\right]$, we obtain (2). The boundedness property (1) is a consequence of the uniform boundedness principle (\cite{PS}, Theorem 2.1.20).

Conversely, assume (1) and (2). Let $f=[f_i;i\in J]\in c_0(J,K)$. Let us prove that $A^{(l)}f\to 0$. For any $\varepsilon >0$, there exists such a finite set $J_0\subset J$ that $|f_j|<\dfrac{\varepsilon}C$ for $j\in J\setminus J_0$. Therefore
\begin{equation}
\left| \sum\limits_{j\in J\setminus J_0}a_{i,j}^{(l)}f_j\right| <\varepsilon \quad \text{for all $i,l$}
\end{equation}
(see Theorem 2.5.1 in \cite{PS} regarding operations with not necessarily countable sequences in non-Archimedean spaces).

On the other hand, let $F=\sup\limits_{j\in J_0}|f_j|$. By virtue of (2), there exists such $l_0\in \mathcal L$ that
$$
\sup\limits_i\left| a_{i,j}^{(l)}\right| <\frac{\varepsilon}F\quad \text{for $l\succ l_0$},
$$
so that
\begin{equation}
\sup\limits_i\left| \sum\limits_{j\in J_0}a_{i,j}^{(l)}f_j\right| <\varepsilon \quad \text{for $l\succ l_0$}.
\end{equation}

It follows from (3) and (4) that $\| A^{(l)}f\| <\varepsilon$ for $l\succ l_0$. $\qquad \blacksquare$

\medskip
The importance of strong operator topology for the theory of operator algebras is explained by the fact that for any set $\mathbf P\subset L(\B )$ of bounded linear operators, the commutant $\mathbf P'=\{ A\in L(\B ):\ AP=PA\text{ for all $P\in \mathbf P$}\}$ is closed in strong operator topology.

\medskip
\section{Group Algebras}

Let $G$ be an arbitrary group. We will consider operators on the Banach space $\B =c_0(G,K)$ of $K$-valued sequences $[x_a;a\in G]$ indexed by elements of $G$. Specifically, we consider the operators $U_g$ of the right regular representation,
$$
U_g[x_a;a\in G]=[x_{ag};a\in G],
$$
and those of the left regular representation
$$
V_g[x_a;a\in G]=[x_{g^{-1}a};a\in G].
$$
These operators are isometric; if
$$
W[x_a;a\in G]=[x_{a^{-1}};a\in G],
$$
then
\begin{equation}
W=W^{-1},\quad WU_gW=V_g.
\end{equation}
As before, with any operator $A\in L(\B )$ we associate a matrix $(\alpha_{a,b})_{a,b\in G}$, infinite, if such is the group $G$.

Denote by $\RR$ the set of all the operators $U_g$, $g\in G$, and by $\LL$ the set of all the operators $V_g$, $g\in G$.

\medskip
\begin{lem}
An operator $A\in L(\B )$ belongs to the commutant $\RR'$, if and only if $\alpha_{a,b}=\eta_{ab^{-1}}$ for all $a,b\in G$, where $\eta:\ G\to K$ is some function. Similarly, $A$ belongs to $\LL'$, if and only if $\alpha_{a,b}=\zeta_{a^{-1}b}$, $a,b\in G$, with some function $\zeta:\ G\to K$.
\end{lem}

\medskip
{\it Proof}. We check directly that
$$
U_g^{-1}AU_g=(\alpha_{ag^{-1},bg^{-1}}).
$$
The fact that $A\in \RR'$ means that
\begin{equation}
\alpha_{a,b}=\alpha_{ag^{-1},bg^{-1}}\quad \text{for any $g\in G$}.
\end{equation}

Set $\eta_c=\alpha_{c,1}$. Choosing $g=b$ in (6) we find that
\begin{equation}
\alpha_{a,b}=\eta_{ab^{-1}}.
\end{equation}
Conversely, the identity (7) implies (6), which proves the assertion regarding $\RR'$.

The automorphism $A\mapsto WAW$ transforms $\RR$ onto $\LL$, and $\RR'$ onto $\LL'$ (see (5)). This leads to the desired description of $\LL'.\qquad \blacksquare$

\medskip
Let $\M ,\N$ be the closed linear spans of $\RR$ (resp. $\LL$) in the strong operator topology.

\medskip
\begin{teo}
\begin{description}
\item[(i)] $\M$ coincides with the closed linear span of $\RR$ in the uniform operator topology. In addition, $\M=\LL'$. Similarly, $\N =\RR'$, and also $\M =\RR''=\M''$, $\N =\LL'' =\N''$.
\item[(ii)] The reduction $\widetilde{\M}$ is isomorphic to the group algebra $\K [G]$.
\item[(iii)] $\M$ is a factor (that is, its center consists of scalar operators), if and only if for each $a\in G$, $a\ne 1$, the set
    $$
    C_a=\left\{ c^{-1}ac,c\in G\right\}
    $$
    of the elements conjugate to a is infinite.
\end{description}
\end{teo}

\medskip
{\it Proof}. (i) Obviously, $\RR \subset \LL'$, so that $\M\subset \LL'$ (since $\LL'$ is closed in the strong operator topology). It suffices to prove that the linear span of $\RR$ is dense in $\LL'$ in the uniform operator topology.

Consider in $c_0(G,K)$ the orthonormal basis $\{ \delta_a,a\in G\}$ (see Section 2), that is $\delta_a=[\delta_{b,a},b\in G]$. We have $U_g(\delta_a)=[\delta_{bg,a},b\in G]=\delta_{ag^{-1}}$. This means that the matrix $(u_{a,b})$ of the operator $U_g$ consists of all zeroes except the ``diagonal'' entries with $a=bg$, which are equal to 1.

On the other hand, by the above Lemma, for the matrix of an element of $\LL'$, there is a constant on every such diagonal. Therefore for $A\in \LL'$ there exists a finite linear combination of operators $U_g$ with a finite set of diagonals coinciding with those of $A$.

Let $\varepsilon >0$. Since $A=(\alpha_{a,b})$ is a bounded operator, for each $b\in G$, there is a finite subset $\mathcal A(b)\subset G$, such that $|\alpha_{a,b}|<\varepsilon$ for $a\notin \mathcal A(b)$. In particular, for $\lambda \notin \mathcal A(1)$, the element $\alpha_{a,b}$ with $a=b\lambda$ coincides with $\alpha_{\lambda ,1}$, and $|\alpha_{\lambda ,1}|<\varepsilon$.

Denote by $A_\varepsilon$ the matrix, in which the elements of the diagonals $a=b\lambda$ with $\lambda \in \mathcal A(1)$ are the same as in $A$, while all other entries are replaced by zero. Then $A_\varepsilon$ is a finite linear combination of operators from $\RR$, and all elements of the matrix $A-A_\varepsilon$ have absolute values $<\varepsilon$. This proves that the linear span of $\RR$ is dense in $\LL'$ in the uniform operator topology; in particular, $\M=\LL'$.

Similarly, $\RR'=\N$, $\N'=\LL'$, and $\N'=\RR''$, so that $\M =\RR'' =\M''$ and $\N =\N''$.

\medskip
(ii) Using the above Lemma, we can identify the algebra $\M$ with the set of sequences $\eta =[\eta_d; d\in G]$ tending to zero by the filter of complements to finite sets, with the $\sup$-norm, the natural structure of a $K$-vector space and the following multiplication operation. Let $\eta ,\zeta \in \M$. Multiplying the matrices $(\eta_{a^{-1}b})_{a,b\in G}$ and $(\zeta_{b^{-1}c})_{b,c\in G}$, we have
$$
(\eta \zeta )_{a,c}=\sum\limits_{b\in G} \eta_{a^{-1}b}\zeta_{b^{-1}c}=\sum\limits_{l\in G}\eta_l\zeta_{l^{-1}(a^{-1}c)},
$$
so that under the above identification, the multiplication of matrices corresponds to the operation
\begin{equation}
(\eta \zeta)_d=\sum\limits_{l\in G}\eta_l\zeta_{l^{-1}d}.
\end{equation}

In the reduction procedure, the unit ball of $\M$ is mapped to the set of $\K$-valued sequences $\{ \widetilde{\eta}_l,l\in G\}$ with finite number of nonzero elements, with the component-wise $\K$-vector space operations and the multiplication induced by (8):
$$
\left(\widetilde{\eta}\widetilde{\zeta}\right)_d=\sum\limits_{l\in G}\widetilde{\eta}_l\widetilde{\zeta}_{l^{-1}d},
$$
that is $\widetilde{\M}$ is isomorphic to $\K [G]$ \cite{P}.

\medskip
(iii) Let $A=(\alpha_{a,b})\in \M \cap \M'$. Then $\alpha_{a,b}=\eta_{a^{-1}b}$ and simultaneously $\alpha_{a,b}=\theta_{ab^{-1}}$ where $\eta ,\theta$ are some $G$-indexed sequences. Taking $b=1$ we find that $\eta_{a^{-1}}=\theta_a$, so that
\begin{equation}
\alpha_{a,b}=\eta_{a^{-1}b}=\eta_{ba^{-1}}.
\end{equation}

Substitute in the second equality in (9) $(b,ab)$ for $(a,b)$. We find that
\begin{equation}
\eta_{b^{-1}ab}=\eta_a,
\end{equation}
so that the sequence $\eta$ is constant on each set $C_a$.

Conversely, this constancy implies that $A=(\eta_{a^{-1}b})\in \M \cap \M'$. To show this, we substitute $ba^{-1}$ for $a$ in the identity (10) and find that $\eta_{a^{-1}b}=\eta_{ba^{-1}}=\theta_{ab^{-1}}$ where $\theta_a=\eta_{a^{-1}}$.

Suppose that $C_a$ is infinite for every $a\ne 1$. We have $\alpha_{1,b}=\eta_b$. Since $A$ is a bounded operator, $\eta_b\to 0$ by the filter of complements to finite sets. This means that $\eta_d=0$ for all $d\in C_a$ ($a\ne 1$). By (10), $\eta_a=0$ for all $a\ne 1$, so that
$$
\alpha_{a,b}=\eta_{a^{-1}b}=\begin{cases}
\eta_1, & \text{for $a=b$},\\
0, & \text{for $a\ne b$},\end{cases}
$$
and $A=\eta_1I$.

Conversely, suppose that $C_{a_0}$ is finite for some $a_0\in G$, $a_0\ne 1$. Define an operator $A$ setting
$$
\alpha_{a,b}=\eta_{a^{-1}b},\quad \eta_d=\begin{cases}
1, & \text{for $d\in C_{a_0}$},\\
0, & \text{for $d\notin C_{a_0}$}.\end{cases}
$$
Then $A=\sum\limits_{d\in C_{a_0}}U_d$ (the finite sum!) and $A\in \M \cap \M'$.

If $\M$ is a factor, then $A=\alpha I$, so that $\alpha_{a,b}=\alpha \delta_{a,b}$. In particular, $\alpha_{1,a_0}=0$. On the other hand, it follows from the definition of $A$ that $\alpha_{1,a_0}=1$. This contradiction shows that $\M$ is not a factor. $\qquad \blacksquare$.

\section{Examples}

In this section we assume that $K=\mathbb C_p$, so that $\K =\mathcal F_p$. We list some classes of groups $G$, for which $\F [G]$ or $Q(\F [G])$ are Baer rings, so that the corresponding operator algebras $\M =\M (G)$ are algebras with Baer reductions.

\medskip
1) {\it Finite groups}. Let $G$ be a finite group. Suppose that $p$ does not divide the order of $G$. It is a consequence of Maschke's and Wedderburn's theorems (see Chapter 3 in \cite {Pi}) that $\F [G]$ is a finite direct sum of full matrix algebras over $\F$. Therefore $\F [G]$ is a type I Baer algebra \cite{Kap}, so that $\M (G)$ is a type I operator algebra with Baer reduction.

\medskip
2) Let $G$ be {\it a torsion free polycyclic-by-finite group}. Then $\F [G]$ has no zero divisors \cite{Cl} and no idempotents different from 0 and 1. Therefore $\F [G]$ is a Baer ring, and $\M (G)$ is a type I operator algebra with Baer reduction.

\medskip
In the following three examples, $\M (G)$ is an operator algebra with the type III weak Baer reduction.

\medskip
3) Let $G$ {\it be locally finite, prime, have no elements of order $p$, and all its conjugacy classes be countable}. Then the ring $Q(\F [G])$ is simple and directly infinite \cite{HO}.

Note that for a locally finite group $G$ without order $p$ elements, the ring $\F [G]$ is regular (see Theorem 1.5 in \cite{P}). In this case, $Q(\F [G])$ is regular and right self-injective (\cite{F}, page 69); every such ring is a Baer ring \cite{G}. For the present example, it follows from Corollary 9.21 in \cite{G} that $Q(\F [G])$ is of type III.

\medskip
4) If $G$ is {\it a solvable group containing a nontrivial locally finite normal subgroup, but no nontrivial finite normal subgroup}, then $Q(\F [G])$ is a type III Baer ring. This result is a combination of Theorem 2.3 of \cite{H79}, Proposition 7 from \cite {H77}, and Theorem 2.10 from \cite{P}.

\medskip
5) Let $G$ be {\it the group of permutations of an infinite set $X$ moving only finitely many elements of $X$}. Then $Q(\F [G])$ is a type III Baer ring; see Corollary 2.8 in \cite{O}.

\medskip
The problem of constructing non-Archimedean operator algebras with type II Baer reductions remains open.

\end{document}